\documentclass[11pt]{article} 
\usepackage{times}
\usepackage{soul}
\usepackage[utf8]{inputenc}
\usepackage{graphicx}
\usepackage{amsmath}
\usepackage{amsthm}
\usepackage{color}
\usepackage{fullpage}
\usepackage{amssymb}
\newcommand{\VSP}{\rule{0ex}{4ex}}
\newcommand\boldpara[1]{\par\noindent{\bf #1.}}
\newcommand\nboldpara[2]{\refstepcounter{#1}\medskip\par\noindent{\large \bf \arabic{#1}. #2.}} 
\newcommand\boldparan[2]{\refstepcounter{#1}\par\noindent{\bf #2 \arabic{#1}.}}
\newcommand\QED{\hfill\rule{1ex}{2ex}}
\newcommand{\tl}{{\symbol{'176}}}
\newcommand\mc{\multicolumn}
\newcommand{\CMV}{\sfsl{cost}_{\sfsls{Mv}}}
\newcommand\sfsl[1]{\textsl{\slshape#1}}
\newcommand\sfsls[1]{\textsl{\scriptsize\slshape#1}}
\newcommand\rms[1]{\textrm{\scriptsize#1}}
\newcommand\bbg[1]{{\mbox{\boldmath{$\bf #1$}}}}
\newcommand\bbgs[1]{{\mbox{\scriptsize\boldmath{$\bf #1$}}}}
\newcommand{\one}{\mathbf{         1}}
\newcommand\defeq{\stackrel{\scriptscriptstyle\rm def}{=}}
\newcommand\myeqno[1]{\refstepcounter{equation}\label{#1}\eqno{(\theequation)}}

\newcommand\optpar{\par}

\newcommand\prm{^T}
\newcommand{\DS}{\displaystyle}

\newcommand{\apile}[2][c]{\begin{array}{@{}#1@{}}#2\end{array}}

\newcommand{\psu}{^{\dagger}}
\newenvironment{shortenum}{\begin{enumerate}\def\theenumii{\roman{enumii}}
\itemsep=0pt \topsep=0pt \parskip=0pt}{\end{enumerate}}

\newcommand{\frc}[2]{{}^{#1\!\!}/_{\!#2}}

\newcommand{\DEL}[1]{}

\begin{document}
\title{On Fast Computation of Directed Graph Laplacian Pseudo-Inverse}

\author{
    Daniel Boley%
    \\
    {\sl University of Minnesota}%
    \\
    {\tt boley@umn.edu}%
}

\maketitle 

\begin{abstract}
The Laplacian matrix and its pseudo-inverse for a
strongly connected directed graph is fundamental
in computing many properties of a directed graph.
Examples include random-walk centrality and betweenness
measures, average hitting and commute times, and other
connectivity measures.
These measures arise in the analysis of many social and
computer networks.
In this short paper,
we show how a linear system involving the Laplacian
may be solved in time linear in the number of edges,
times a factor depending on the separability of the graph.
This leads directly to the column-by-column computation of the
entire Laplacian pseudo-inverse
in time quadratic in the number of nodes, i.e., constant time per matrix entry.
The approach is based on ``off-the-shelf'' iterative methods
for which global linear convergence is guaranteed, without
recourse to any matrix elimination algorithm.
\end{abstract}

\boldpara{Keywords:} 
Graph Laplacian;
Directed Graphs;
Pseudo-Inverse;
Iterative Methods.


\nboldpara{par}{Introduction} 
\optpar
Many properties of networks 
can be found via the solution of special linear systems
based on the graph Laplacian.  Examples 
include the well-known pagerank, centrality measures, betweenness measures, graph cuts,
distances or affinities between nodes, trust/influence propagation,
etc.\ \cite{%
boley11,%
ChungDigraph,%
Fouss07,%
golnari2014revisiting,
Liu16,%
ZhouHuangScholkopf05%
}.  
These properties have spawned many papers on efficient,
almost linear time solvers for these special linear systems
such as \cite{ST14} for symmetric systems for undirected graphs to
more recent papers reporting almost linear time for
non-symmetric Eulerian Laplacians for directed graphs
\cite{Cohen18,Cohen16}. 
{For the purposes of this paper, we say a Laplacian matrix $\mathbf L$
is ``Eulerian'' if $\mathbf L$ has nullity 1 and $\mathbf L\mathbf w = \mathbf L^T\mathbf w = \mathbf 0$ for a strictly
positive vector $\mathbf w$.}
These fast methods use a careful
ordering of the nodes, an approximate factorization using
Gaussian elimination used as a preconditioner to an
iterative method based on, e.g., Richardson iteration.
The theoretical
running time for the methods of 
\cite{Cohen18,Cohen16}
can be bounded by $O(m)\log^c(n\kappa\varepsilon)$
with high probability $(1-\delta)$, where $m$ is the number of edges, $n$ the number of
nodes, $\kappa$ is the matrix condition number, and $\varepsilon$ is 
the desired accuracy, with $O(n)\varepsilon^{-2}\log^c(1/\delta)\log^c(n\kappa\varepsilon)$
fill-in from the
inexact factorization.  The $c$'s are some arbitrary constants.
In this short paper we use a different approach to obtain an algorithm for the pseudo-inverse of
a non-symmetric Eulerian Laplacian.  Our approach is to use only iterative methods in widespread use in practice,
and which also enjoy provable linear
convergence guarantees and per-iteration costs linear in the number of edges
in the graph.  We also propose a computational process to obtain an Eulerian scaling.
By using only iterative methods, we avoid the issue of fill-in entirely.
{This paper focuses on strongly connected directed graphs.
The results carry over to the case of connected undirected graphs,
but most of the results can be simplified.  This is beyond the scope of this paper.}

The rest of this paper begins with
preliminaries to introduce the Laplacians and other
basic matrices associated with directed graphs,
followed by 
a {theorem} 
which reduces the pseudo-inverse computation to a simple matrix inversion.
Then we present the overall algorithm to find an Eulerian scaling and compute the
pseudo-inverses for Eulerian Laplacians, followed by an outline of the complexity analysis, which includes the
convergence theory and the cost per iteration.
We then briefly show how the pseudo-inverse for a non-Eulerian
Laplacian can be recovered from that of an Eulerian Laplacian.
We end with a short table of experiments showing the performance of
the methods in practice is consistent with the theoretical complexity bounds.
We collect existing theoretical results on which our methods
are based into an Appendix.
Throughout this paper, all norms are the matrix or vector
2-norms, unless otherwise specified.

\nboldpara{par}{Preliminaries} 
Consider a directed graph with adjacency matrix $\mathbf A \in \mathbb R^{n \times n}$
where $a_{ij}$ is the weight on the edge $i ${$\rightarrow$}$ j$ if such an edge
exists, otherwise $a_{ij}= 0$.  If $\one $ is the vector of all ones of
appropriate dimension, then $\mathbf d = \mathbf A \one$ is the vector of out-degrees, 
$\mathbf{D}= \mathsf{Diag}(\mathbf d)$ is the diagonal matrix with the entries of $\mathbf d$
on the diagonal, and
$\mathbf{P}= \mathbf{D}^{-1} \mathbf A$ is the matrix of transition
probabilities for a random walk over this directed graph.
Throughout this paper we assume the graph is strongly connected implying that $\mathbf P$ is irreducible.
Let $\bbg\pi$ be the unique vector of stationary probabilities over this graph,
i.e., the vector satisfying $\bbg\pi\prm \mathbf P = \bbg\pi\prm$ and $\bbg\pi\prm\one  = \one $,
and let $\bbg \Pi = \mathsf{Diag}(\bbg\pi)$ be the diagonal matrix with the
stationary probabilities $\{\pi_i\}_1^n$ on the diagonal.
Perron-Frobenius theory guarantees $\bbg \pi$ exists and is strictly positive
\cite{GantmacherII,hojo:91}.
{Several different Laplacians have been defined for a given digraph, each related 
to each other through a variety of diagonal scalings \cite{boley11}:}
\begin{equation} \label{unsymlap}
\apile[llll]{
\mathbf L^{\mathrm r} & = & \bbg \Pi - \bbg \Pi \mathbf P & \mbox{\small random walk Laplacian} \\
\mathbf L^{\mathrm a} & = & \mathbf D - \mathbf A = \mathbf D - \mathbf D \mathbf P & \mbox{\small unnormalized Laplacian}\\
\mathbf L^{\mathrm p} & = & \mathbf I - \mathbf P & \mbox{\small normalized Laplacian} \\
\mathbf L^{\mathrm d} & = & \mathbf I - \bbg \Pi^{\frc12} \mathbf P \bbg \Pi^{-\frc12} & \mbox{\small diagonally scaled
Laplacian}
}
\end{equation}
and corresponding pseudo-inverses
\begin{equation} \label{psunsymlap}
\mathbf M^{\mathrm r} = (\mathbf L^{\mathrm r})\psu,
~ \mathbf M^{\mathrm d} = (\mathbf L^{\mathrm d})\psu ,
~ \mathbf M^{\mathrm p} = (\mathbf L^{\mathrm p})\psu,
~ \mbox{etc.}
\end{equation}
{It is well known that the Laplacians are interchangeable in the limited sense that
one can obtain many graph properties from
one or another of the
Laplacians or their pseudo-inverses,}
e.g.,
the average length $h(i,k)$ of a random walk starting from node $i$ before reaching node $k$
and the average round-trip commute time $c(i,k)$ 
\cite{Aldous14,
boley11,
Fouss07,
GrinsteadSnell,
KemenySnell76,
LiZhang10a,
LiZhang10b,
Norris97}
(even for strongly connected directed graphs):
$${
h(i,k) 
{~ =  \mathbf m_{kk}^{\mathrm r}-\mathbf m_{ik}^{\mathrm r} + \sum_\ell (\mathbf m_{i\ell}^{\mathrm r}-\mathbf m_{k\ell}^{\mathrm r})\pi_\ell}
 = \frac{\mathbf m_{kk}^{\mathrm d}}{\pi_k} 
 - \frac{\mathbf m_{ik}^{\mathrm d}}{\sqrt{\pi_i\pi_k}}
}\myeqno{transit-times-m}
$$
$${
c(i,k) 
{~
 = \mathbf m_{kk}^{\mathrm r}+\mathbf m_{ii}^{\mathrm r}-\mathbf m_{ik}^{\mathrm r}- \mathbf m_{ki}^{\mathrm r}}
 = \frac{\mathbf m_{kk}^{\mathrm d}}{\pi_k}
  +\frac{\mathbf m_{ii}^{\mathrm d}}{\pi_i}
  -\frac{\mathbf m_{ik}^{\mathrm d}+ \mathbf m_{ki}^{\mathrm d}}{\sqrt{\pi_i\pi_k}}
}
$$
{The choice of which Laplacian scaling to use
depends on which leads to a simpler formula.} 

The pseudo-inverse also yields the average number of visits to an individual node $j$
for random walks starting in node $i$ before reaching $k$:
\cite{boleyrandom,LAA-2019}:
$$
v(i,j,k)
=  (\mathbf m_{ij}^{\mathrm r} - \mathbf m_{kj}^{\mathrm r} - \mathbf m_{ik}^{\mathrm r} + \mathbf m_{kk}^{\mathrm r}) \pi_j
{= \sqrt{\frac{\pi_j}{\pi_i}} \mathbf m_{ij}^{\mathrm d}
  - \sqrt{\frac{\pi_j}{\pi_k}} \mathbf m_{kj}^{\mathrm d} 
  - \frac{\pi_j}{\sqrt{\pi_i\pi_k}} \mathbf m_{ik}^{\mathrm d}
  + \frac{\pi_j}{\pi_k}\mathbf m_{kk}^{\mathrm d}}
,
\myeqno{visits-m}$$
and the probability that such a random walk passes node $j$ at all:
$$\mathsf{Prob}(\mbox{pass }j \mbox{ on walks } i \mbox{$\rightarrow$} k) =
v(i,j,k)/v(j,j,k) . \myeqno{prob-visit}$$
By summing $v(i,j,k)$ across various dimensions, one can obtain various centrality and betweenness measures
for individual nodes \cite{LAA-2019}.  For example, it can be easily verified that summing (\ref{visits-m})
over $j$ yields formula (\ref{transit-times-m}):
$\sum_j v(i,j,k) = h(i,k)$, and summing (\ref{transit-times-m}) over $k$ yields
$\sum_k \pi_k h(i,k) = \textsf{trace } \mathbf M^{\mathrm d}$ (independent of $i$)
equivalent to the Random Target Lemma \cite{Aldous14}.

\boldpara{Remark}
To illustrate how these relations are a simple consequence of our Theorem \ref{psu} below,
a derivation for the last equality in (\ref{visits-m}) is given here.  A derivation
had not given previously elsewhere to the author's knowledge.
Apply 
{Theorem} \ref{psu} by setting
$A = \mathbf L^{\mathrm d} = \bbg \Pi^{\frc12} (\mathbf I -  \mathbf P) \bbg \Pi^{-\frc12} $,
$\mathbf u = \sqrt{\bbg \pi}$, and
$B = \mathbf M^{\mathrm d} = (\mathbf L^{\mathrm d} )\psu$.
Writing ($\ref{AA})$ in Theorem \ref{psu}b elementwise, we have:
$${
[A_{11}^{-1}]_{ij} 
= \mathbf m_{ij}^{\mathrm d}
  - \sqrt{\frac{\pi_i}{\pi_n}} \mathbf m_{nj}^{\mathrm d} 
  - \sqrt{\frac{\pi_j}{\pi_n}} \mathbf m_{in}^{\mathrm d}
  + \frac{\sqrt{\pi_i\pi_j}}{\pi_n}\mathbf m_{nn}^{\mathrm d}
}$$
Apply two-sided diagonal scaling 
$(\mathbf I -  \mathbf P) = \bbg \Pi^{-\frc12} A \bbg \Pi^{\frc12} $ to obtain
$${ v(i,j,n) = [(\mathbf I - \mathbf P_{11})^{-1}]_{ij} =
\sqrt{\frac{\pi_j}{\pi_i}}[A_{11}^{-1}]_{ij} 
= \sqrt{\frac{\pi_j}{\pi_i}} \mathbf m_{ij}^{\mathrm d}
  - \sqrt{\frac{\pi_j}{\pi_n}} \mathbf m_{nj}^{\mathrm d} 
  - \frac{\pi_j}{\sqrt{\pi_i\pi_n}} \mathbf m_{in}^{\mathrm d}
  + \frac{\pi_j}{\pi_n}\mathbf m_{nn}^{\mathrm d}
}$$
\QED

Given an arbitrary directed graph with $n$ nodes, one can augment the graph with an extra node $n${$+$}$1$ 
such that, upon every transition in a random walk over the graph, there is a small probability $\gamma$ that the walker
"evaporates" to node $n${$+$}$1$, and thence the walker transitions to an arbitrary node with equal probability
(or biased probabilities in a personalized setting).
This is a process very similar to teleportation in the pagerank setting.
The result is a strongly connected directed graph to which we can apply the methods of this paper.
In this case the average number of visits $v(i,j,n\mbox{$+$}1)$ or average path lengths $h(i,j)$
would yield a continuum of affinity estimates from 
$i$ to $j$, approximating random walk affinity for $\gamma$ near 0 and shortest path affinity for $\gamma$ near 1.
{The probability $v(i,j,n\mbox{$+$}1)/v(j,j,n\mbox{$+$}1)$ of equation
(\ref{prob-visit}) can be interpreted as the trust of node $j$ from the point
of view of node $i$ in a trust network \cite{Buendia17,Liu16}.  The sum 
$\sum_i v(i,j,n\mbox{$+$}1)/v(j,j,n\mbox{$+$}1)$ can be interpreted as
an average level of trust in node $j$ or a
measure of influence of node $j$ in a social network
\cite{golnari2014revisiting}.}

The main contributions of this paper are: 
(A) we show how an off-the-shelf iterative method in widespread use yields a method to find
the pseudo-inverse of
an Eulerian Laplacian with a provable complexity guarantee that is linear in the number of edges times a factor related to
the connectness of the graph ; (B) we show how another off-the-shelf method yields a method to find a Eulerian scaling
for a non-Eulerian Laplacian with similar complexity guarantees;
and (C) we illustrate the methods with some examples showing
the linear complexity can be observed in practice with the off-the-shelf numerical procedures.
The constructions proposed in this paper are kept as simple as possible
to highlight a minimal set of assumptions needed to form the basis for a fast Laplacian solver.
Most of the theoretical properties used in this paper are well-known, but 
we include a few brief proofs to make this paper more self-contained.

{The computation of the Moore-Penrose pseudo-inverse in the general case
usually
proceeds using the SVD using $O(n^3)$ time \cite{GvL13}.  However, a graph
Laplacian for strongly connected digraphs has nullity equal to 1. This
special property allows one to compute the pseudo-inverse with just an ordinary
matrix inversion using the formulas given in theorems \ref{psu}(a),
\ref{non-euler}(a) 
below,
based on the theory of \cite{Meyer73},
assuming one knows the left and right nullspaces for the Laplacian.
The cost of the matrix inversion using a standard algorithm like
Gaussian elimination is also 
$O(n^3)$ but much faster than a full SVD \cite{GvL13}.
{The cost to obtain even a single column of the pseudo-inverse
using Gaussian elimination is still $O(n^3)$.}
The goal
in this paper is to show how this complexity can be reduced 
to $O(m\cdot\log\kappa\varepsilon)$ for each
column of the pseudo-inverse and to $O(nm\cdot\log\kappa\varepsilon)$ for the entire pseudo-inverse,
where $m$ is the number of edges.
The procedures outlined here also include the computation of the necessary left and
right nullspaces with $O(m\cdot\log\kappa\varepsilon)$ cost.
{An undirected graph would lead to subtantial simplifications and often
lower cost using different techniques, but this is beyond the scope of this paper.} 
} 

\nboldpara{par}{{Theoretical} Construction}
\optpar
In this paper we study mainly the Eulerian Eulerian Laplacian matrices $\mathbf L^{\mathrm r}$ and $\mathbf L^{\mathrm d}$
(\ref{unsymlap}).
The matrix $\mathbf L^{\mathrm r}$ can be thought of as the unnormalized Laplacian for a weighted digraph with
adjacency matrix $\bbg \Pi \mathbf P$.  This last matrix has all row sums and column sums
equal to each other: $\bbg \Pi \mathbf P \one  = ( \bbg \Pi \mathbf P)^T  \one  = \bbg \pi$, 
and hence the corresponding Laplacian is ``Eulerian'' \cite{Cohen18,Cohen16}.
The matrix $\bbg \Pi^{\frc12} \mathbf P \bbg \Pi^{-\frc12}$ has a similar property:
$\bbg \Pi^{\frc12} \mathbf P \bbg \Pi^{-\frc12} \sqrt{\bbg\pi}$ $=$
$(\bbg \Pi^{\frc12} \mathbf P \bbg \Pi^{-\frc12})^T \sqrt{\bbg\pi} $ 
$=$ $\sqrt{\bbg\pi} $.
Here 
$\sqrt{\bbg\pi} $ $=$ $[\sqrt{\pi_1},\ldots,\sqrt{\pi_n}]^T$.

{In this paper, we focus specifically on Laplacian matrices corresponding to
strongly connected digraphs, specifically
matrices $L$ such that are {\sf (Pa)} irreducible, {\sf (Pb)} have all positive diagonal
entries and no positive off-diagonal entries, and {\sf (Pc)} satisfy
$L \mathbf x = 0$ for some strictly positive vector $\mathbf x > 0$.  We call such a matrix  an
{\em Eulerian Laplacian} if $L^T \mathbf x = L \mathbf x = 0$ for some strictly positive vector $\mathbf x > 0$.  }

The main point of this section is to present the mapping between the pseudo-inverse computation
for an Eulerian Laplacian
and the computation of related ordinary inverses.  In the following {theorem} we present two such mappings.
Part (a) connects the pseudo-inverse of the Eulerian Laplacian matrix with the ordinary inverse of a
symmetric rank-1 modification to that matrix.
{This is a special case of the general theory from
\cite{Meyer73}.}
The rank-1 modification is exactly in the direction corresponding
to the nullspace of the original Laplacian.  This construction is well known
(see, e.g., \cite{boley11,Fouss07}).
Part (b) shows how the ordinary inverse of the $(n-1)\times(n-1)$ principal submatrix of a Eulerian Laplacian
can be obtained directly from the pseudo-inverse of the entire matrix via small rank changes, and part (c)
gives reverse mapping, from the ordinary inverse of the submatrix to the pseudo-inverse of the entire matrix.
These connections will allow the use of off-the-shelf iterative methods for the ordinary inverse in order to 
obtain the desired pseudo-inverse.
{Part (b) is a special case of \cite[Lemma 1]{boley11}. We remark that for an Eulerian Laplacian, the Moore-Penrose pseudo-inverse
is the same as the group inverse \cite{Meyer75,Roberts-GroupInv},
and hence part (c) of {Theorem} \ref{psu} is a special case
of \cite[Thm 5.2]{Meyer75}.}

\boldparan{thm}{{Theorem}} \label{psu}
Let $C$ be an $n \times n$ non-singular matrix and suppose
$ A = C - \alpha \mathbf u \mathbf u\prm$
is singular with $A\mathbf u = A^T \mathbf u = \mathbf 0$. 
Partition $A$ and $\mathbf u$  as follows:
$${
A = \left[\apile[cc]{A_{11} & \mathbf a_{12} \\ \mathbf a_{21}^T & a_{nn}}\right]
\mbox{,~~~}
\mathbf u = \left[\apile[c]{ \mathbf u_{1} \\ u_{n}}\right] .
}\myeqno{partition}$$
Assume $\mathbf u\prm \mathbf u = 1$ and $u_n > 0$.
Then
\begin{itemize}
\item[(a)]
the left and right nullspaces of $A$ are 
$\mathsf{nullsp}(A)=\mathsf{nullsp}(A^T)=\mathsf{span}(\mathbf u)$, and
the Moore-Penrose pseudo-inverse of $A$ is given as:
\vspace{-1ex}
$${
A\psu ~\defeq~ B ~ = {C^{-1}}- \mathbf u \mathbf u^T / \alpha.
}$$
\vspace{-3ex}
\item[(b)]
$A_{11}^{-1}$ exists and can be written in terms of $A\psu = B$:
$${
\apile[@{}l     l     l@{}]{
A_{11}^{-1}
&=& 
\left[ I_{n-1} ,~ -\mathbf u_1/u_n \right]
B
\left[\apile{ I_{n-1} \\ -\mathbf u_1^T/u_n} \right]
\\[3ex] &=& 
   B_{11} - \frac1{u_n}{\mathbf u_1 } \mathbf b_{21}^T - \frac1{u_n}\mathbf b_{12} {\mathbf u_1^T}
    + \frac{b_{nn}}{u_n^2} \mathbf u_1 \mathbf u_1^T \mbox{~~~} ,
}
}\myeqno{AA}$$
where we have partitioned $B$ as in (\ref{partition}).
\item[(c)] We can write $A$ and $B = A\psu$ in terms of $A_{11}$ and $\mathbf u$
as follows
$${
\apile[@{}*{4}{l     }l@{}]{
A 
&=&
{\left[ \apile[cc]{A_{11}  & -\frac1{u_n} A_{11} {\mathbf u_1}
       \\ -\frac1{u_n}{\mathbf u_1^T} A_{11} & \frac1{u_n^2}{\mathbf u_1^T A_{11} \mathbf u_1}
       }\right]}
&=& 
{\left[\apile{ I_{n-1} \\ -\frac1{u_n}{\mathbf u_1^T}} \right]}
A_{11}
{\left[ I_{n-1} ,~ -\frac1{u_n}{\mathbf u_1} \right] ;}
\\[-.5ex]
\\[-.5ex]
B
&=&
  \left[\apile[cc]{
  B_{11} & \mathbf b_{12}
  \\
   \mathbf b_{21}^T
   & 
   b_{nn} 
  } \right]
&=&
\left[\apile{ I_{n-1} \mbox{$-$} \mathbf u_1\mathbf u_1^T \\ -u_n\mbox{$\cdot$}\mathbf u_1^T} \right]
A_{11}^{-1}
\left[ I_{n-1} \mbox{$-$} \mathbf u_1\mathbf u_1^T , -u_n\mbox{$\cdot$}\mathbf u_1 \right] 
}
}\myeqno{BB}$$
where the individual blocks are 
$${
\apile[l@{}l     l     l@{}]{
&B_{11} &=& A_{11}^{-1} - \mathbf u_1 \mathbf t^T 
- \mathbf w \mathbf u_1^T  + (\mathbf u_1^T \mathbf w)\mbox{$\cdot$}\mathbf u_1 \mathbf u_1^T 
\\&
{\mathbf b_{12}} &=& {u_n(\mathbf u_1^T \mathbf w ) \mbox{$\cdot$}\mathbf u_1 
   - u_n\mbox{$\cdot$} \mathbf w}
\\&
\mathbf b_{21}^T &=&
   u_n(\mathbf u_1^T    \mathbf w          ) \mbox{$\cdot$}\mathbf u_1^T - u_n\mbox{$\cdot$}\mathbf t^T
\\&
b_{nn} &=& u_n^2 (\mathbf u_1^T \mathbf w )
\\ \mc{4}{@{}l}{\mbox{where}
~~~ \mathbf w = A_{11}^{-1} \mathbf u_1 ,
~~~ \mathbf t^T = \mathbf u_1^T  A_{11}^{-1} } .
}
}$$

\end{itemize}
\boldpara{Proof}
\begin{itemize}
\item[(a)] A simple calculation yields $AB = BA = I_n - \mathbf u \mathbf u^T$, and a further simple calculation
yields $ABA = A$ and $BAB = B$.  Hence $B$ satisfies the conditions to be the Moore-Penrose pseudo inverse.
\item[(b)] $A \mathbf u = \mathbf 0$ and $\mathbf u^T A = \mathbf 0^T$ imply
$${
\apile[l     l     l]{ 
A_{11} \mathbf u_1 &=& - \mathbf a_{12} u_n, 
\\
\mathbf a_{21}^T \mathbf u_1 &=& -  a_{nn} u_n, 
}
\quad
\quad
\apile[l     l     l]{
\mathbf u_1^T A_{11} &=& - \mathbf a_{21}^T u_n
\\
\mathbf u_1^T \mathbf a_{12} &=& -  a_{nn} u_n
}
}$$
Likewise,
$B \mathbf u = \mathbf 0$ and $\mathbf u^T B = \mathbf 0^T$ imply
$${
\apile[l     l     l]{
B_{11} \mathbf u_1 &=& - \mathbf b_{12} u_n ,
\\
\mathbf b_{21}^T \mathbf u_1 &=& -  b_{nn} u_n ,
}
\quad\quad
\apile[lll]{
\mathbf u_1^T B_{11} &=& - \mathbf b_{21}^T u_n
\\
\mathbf u_1^T \mathbf b_{12} &=& -  b_{nn} u_n
}
}$$
These yield the equivalence for the two formulas for $A_{11}^{-1}$ in (\ref{AA}) 
and the formulas for $A$ in
(\ref{BB}). 
To verify (\ref{AA}) is indeed the inverse of $A_{11}$,
we multiply (\ref{AA}) by $A_{11}$ to obtain the identity:
$${
\apile[l     l     l]{
\mc{3}{l}{\left[ I_{n-1} ,~ -\mathbf u_1/u_n \right]
B
\left[\apile{ I_{n-1} \\ -\mathbf u_1^T/u_n} \right]
\cdot A_{11}} \mbox{~~~~~~~~}
\\ \mbox{~~~~~} &=& 
\left[ I_{n-1} ,~ -\mathbf u_1/u_n \right]
B
\left[\apile{ A_{11}  \\ -\mathbf a_{21}^T} \right]
\\\mbox{~~~~~} &=& 
\left[ I_{n-1} ,~ -\mathbf u_1/u_n \right]
\left[\apile{ I_{n-1} - \mathbf u_1 \mathbf u^T \\ - u_n \mathbf u_{1}^T} \right]
= I_{n-1}
}
}$$
\item[(c)]
Using the second formulas for $A$, $B$ in (\ref{BB}), calculations similar to the proof of (b) yield
$AB = BA = I_n - \mathbf u \mathbf u^T$ and then $ABA = A$, $BAB = B$.
\QED
\end{itemize}

\nboldpara{par}{Algorithm for Eulerian Laplacians}
\optpar
We study the computation of the pseudo-inverses of $\mathbf L^{\mathrm d}$ and $\mathbf L^{\mathrm r}$.
Using {Theorem} \ref{psu}, we can write these as follows:
$${\apile[l     l     l     l     l     l]{
{}[\textsf{a}] & \mathbf M^{\mathrm d} & \defeq & (\mathbf L^{\mathrm d})\psu &=& (\mathbf L^{\mathrm d}
+\sqrt{\bbg \pi}\sqrt{\bbg \pi}^T)^{-1}- \sqrt{\bbg \pi}\sqrt{\bbg \pi}^T 
\\
{}[\textsf{b}] & \mathbf M^{\mathrm r} & \defeq & (\mathbf L^{\mathrm r})\psu &=& \left(\mathbf L^{\mathrm r}
+ \alpha \frac{\one \one ^T}{n}\right)^{-1} -\frac{1}{\alpha} \frac{\one  \one ^T }{n}},
}\myeqno{LapInv}$$
for some arbitrary $\alpha\not=0$ (we use $\alpha=1$ below, but include it here to show a slightly more general formula).
{Theorem} \ref{psu} applies here because
both $\mathbf L^{\mathrm r}$ and $\mathbf L^{\mathrm d}$ are Eulerian.

The overall algorithm begins with a computation of the stationary probabilities.  These probabilities
are used to scale non-Eulerian Laplacians to an Eulerian scaling.  The final step is to solve for the pseudo-inverse
of the Eulerian Laplacian by applying an iterative method to (\ref{LapInv}).  The detailed steps
are given in Algorithm \ref{alg1}.

\boldparan{alg}{\VSP Algorithm} \label{alg1} \\
{\bf Input:} $\mathbf P$ $=$ probability transition matrix for a random walk
over the graph {and an index set ${\cal J} \subset \{1,\ldots,n\}$.}
\\{\bf Output:} Stationary probabilities $\bbg \pi$ and {the columns indexed by
${\cal J}$ of either
{}[\textsf{a}] pseudo-inverse $\mathbf M^{\mathrm d} = (\mathbf L^{\mathrm d} )\psu = (\mathbf I - \bbg \Pi^{\frc12} \mathbf P \bbg \Pi^{-\frc12})\psu$,
or {}[\textsf{b}] pseudo-inverse $\mathbf M^{\mathrm r} = (\mathbf L^{\mathrm r} )\psu = (\bbg \Pi - \bbg \Pi \mathbf P )\psu$.}
Note: items marked {}[\textsf{a}] are needed only for $\mathbf M^{\mathrm d} $ while items marked {}[\textsf{b}] are needed only 
for $\mathbf M^{\mathrm r}$.
\begin{shortenum}

   \item \label{a1step1} Compute $\bbg \pi$, the vector of stationary probabilities:
   \\
	   Use the modified subspace iteration method with $\ell$ starting vectors \cite{Stewart} on $\mathbf P^T$ to compute eigenvector corresponding to
      the eigenvalue $\lambda = 1$.  Here $\ell$ is larger than the period of the graph.


   \item  \label{a1step2}
   Set {}[\textsf{a}] $\mathbf L^{\mathrm d} = \mathbf I - \bbg \Pi^{\frc12} \mathbf P \bbg \Pi^{-\frc12}$,
   where $\bbg \Pi^{\frc12} = \mathsf{Diag}(\sqrt{\bbg \pi}) $,
   or\\ \mbox{~~~~~}
   {}[\textsf{b}] $\mathbf L^{\mathrm r} = \bbg \Pi - \bbg \Pi \mathbf P $,
   where $\bbg \Pi = \mathsf{Diag}(\bbg \pi) $.

   \item Compute {the selected columns indexed by ${\cal J}$ of
   pseudo-inverse of the Eulerian Laplacian using (\ref{LapInv}): \label{a1step5}}
   \\ \mbox{{either }} [\textsf{a}] $\mathbf M^{\mathrm d}  =   (\mathbf L^{\mathrm d})\psu$
   \mbox{or} [\textsf{b}] $\mathbf M^{\mathrm r}  =   (\mathbf L^{\mathrm r})\psu$ column-by-column as follows:
      \\
	For $j \in {\cal J}$:
      {\leftmarginii=7ex%
      \begin{shortenum}
         \item \label{a1inv} Solve the following linear systems using restarted GMRES($\ell$)
	 for $\mathbf x^{\mathrm d}$ and $\mathbf x^{\mathrm r}$:
	 \\ \mbox{~~~~~}[\textsf{a}] $(\mathbf L^{\mathrm d} + \sqrt{\bbg \pi}\sqrt{\bbg \pi}^T) \mathbf x^{\mathrm d} = \mathbf e_j$,
	  or \\ \mbox{~~~~~}[\textsf{b}] $(\mathbf L^{\mathrm r} + \one  \one ^T/n ) \mathbf x^{\mathrm r} = \mathbf e_j$.
         \item \label{a1fill} Fill in the $j$-th column of pseudo-inverse:
	 \\ \mbox{~~~~~}[\textsf{a}]~$\mathbf M_{:,j}^{\mathrm d} = \mathbf x^{\mathrm d} - \sqrt{\pi_j}\cdot \sqrt{\bbg \pi}$,
	 or
	 \\ \mbox{~~~~~}[\textsf{b}] $\mathbf M_{:,j}^{\mathrm r} = \mathbf x^{\mathrm r} - \frc1n \one $,
      \end{shortenum}
	   }
\end{shortenum}

\nboldpara{par}{Complexity of Algorithm \ref{alg1}: Convergence}
\optpar
The two most expensive steps in Algorithm \ref{alg1} are steps \ref{a1step1} and \ref{a1step5}(i), both involving an iterative
method. Their cost is a product of the cost per iteration times the number of iterations.
{Of the remaining steps,
step \ref{a1step2} involves diagonal scaling which costs only $O(m)$
operations, since only the nonzero elements must be computed.
Step \ref{a1step5}(ii) costs $O(n)$ for each column or $O(n|{\cal J}|)$ altogether.
If the entire pseudo-inverse
were to be computed,} it could cost $O(n^2)$ overall, i.e.,
constant time per matrix entry.

In step \ref{a1step1} the modified subspace iteration (see Algorithm \ref{algsub} in the Appendix)
\cite{Stewart} computes the Schur decomposition of a small $\ell \times \ell$ matrix which is the orthogonal projection of the original
matrix $\mathbf P^T$ onto an $\ell$ dimensional subspace. 
If $\lambda_1,\lambda_2,\ldots,\lambda_n$ are the eigenvalues of $\mathbf P$ with
$1 = \lambda_1 \ge |\lambda_2| \ge |\lambda_{\ell}| > |\lambda_{\ell+1}| \ge \cdots \ge |\lambda_n|$,
and $\lambda_1=1$ is a simple eigenvalue,
then Stewart \cite{Stewart} showed that the leading eigenvector
(corresponding to $\lambda_1=1$) converges as:
$${
\frac{\|\mathbf P^T \mathbf x^{[k]} - \bbg x^{[k]}\|} {\|\mathbf P^T \mathbf x^{[0]} - \bbg x^{[0]}\| }
\le O(|\lambda_{\ell+1}|^k) 
}$$
where $\mathbf x^{[k]}$ denotes the approximation to the eigenvector corresponding to $\lambda_1=1$ at the $k$-iteration.
In the following we use $c_1,c_2,\ldots$ to represent small positive constants in the costs bounds,
all of which are less than 10.
If the random walk is aperiodic, then we are guaranteed that $1 > |\lambda_2| \ge \cdots \ge |\lambda_\ell| \ge \cdots $.
If the random walk is periodic with period $\sfsl{per}$, it suffices to have $\ell > \sfsl{per}$
in order to guarantee that $\lambda_{\ell+1} < \lambda_1 =1 $.
To obtain 
an error at most $\sfsl{tol}\,$  requires
at least $c^{\mathrm s}$ iterations with 
$${
c^{\mathrm s} \ge \frac{ |\log(\sfsl{tol})|+|\log ( \|\mathbf P^T \mathbf x^{[0]} - \mathbf x^{[0]}\| )|}{|\log(\lambda_{\ell+1})|} 
}$$
As written in Algorithm \ref{algsub} in the Appendix,
{the cost per iteration is {$\sfsl{cost}_1^{\rms{subspace}} = O(\CMV\cdot(\ell+1) + (n\ell^2)
+ (\ell^3))$}, where $\CMV$ is the cost of one matrix-vector product, proportional to the
number of nonzero entries in the matrix, $\mathsf{nnz}(\sfsl{matrix})$.
The storage required is
$\sfsl{space}^{\rms{subspace}} =  O(\sfsl{\#edges} + n) +  (n \ell) + \ell^2)$
We remark that for undirected graphs, this eigenvector is a multiple of the vector of
degrees, so this step would be essentially free. }

We remark that there are many choices of algorithms to compute this eigenvector, similar to the many choices to compute the
pagerank vector, many of which can be faster \cite{GZB05}. 
{If the dimension is small enough, one can use solve for the eigenvector directly by finding $\mathbf v^T$ satisfying the homogeneous linear system
$${
[\mathbf v^T,1] \left[\apile{\mathbf P_{11} \\ \mathbf p_{21}^T} \right] = \mathbf v^T,
}\myeqno{directeig}$$
where $\mathbf P_{11}$ is the upper left $(n-1) \times (n-1)$ block of $\mathbf P$, but
with $O(n^3)$ cost with $O(n^2)$ space using ordinary Gaussian elimination.}
{Modified Subspace Iteration is an effective algorithm which enjoys a simple
bound on its convergence rate and fixed cost per iteration and
little additional space beyond that of the input matrix, leading to
a simple complexity bound.}

The other costly step is step \ref{a1step5}(i) to compute a column of the inverse $\mathbf M^{\mathrm d}$. 
This line is called $|{\cal J}|$ times, each time solving a linear system to obtain one
column of the inverse. {To solve the linear system, we use an iterative
method with a fixed bound on the cost per iteration and a convergence
guarantee yielding a bound on the number of iterations depending on the
accuracy desired but not on the dimensions of the graph.
For this purpose we use
GMRES({$\ell$}), i.e., restarted GMRES where $\ell$ is the number of inner steps
between restarts, because it not only enjoys these theoretical properties,
but has also been observed to be a very effective solver in practice
\cite{Saad-Iter,Saad+Schultz}.  }
{The cost of one outer iteration of
restarted GMRES is
(details in the Appendix)
$\sfsl{cost}_1^{\rms{GMRES}} = 
O(\ell\cdot\CMV + n\ell^2+\ell^{2}) $.}
{In order to complete the complexity bound for this step,
we must show that restarted GMRES
converges at a guaranteed rate.}
In order to do that, we show that the symmetric part of the modified Laplacian
matrices in question
are positive definite, in the following lemma.

\boldparan{thm}{Lemma}
If $\mathbf P$ is the probability transition matrix for a strongly connected directed graph,
and $\bbg \pi>0$ is the vector of stationary probabilities, then the following two matrices
$${ \apile[l     l     l]{
S(\mathbf L^{\mathrm d}) 
+ \sqrt{\bbg \pi}\sqrt{\bbg \pi}^T
&\defeq & (\mathbf L^{\mathrm d} +  (\mathbf L^{\mathrm d} )^T)/2 
+ \sqrt{\bbg \pi}\sqrt{\bbg \pi}^T
\\
S(\mathbf L^{\mathrm r})
+ \one  \one ^T 
&\defeq & (\mathbf L^{\mathrm r} +  (\mathbf L^{\mathrm r} )^T)/2 
+ \one  \one ^T 
}}$$
are symmetric positive definite.

\boldpara{Proof (sketch)}
We show the symmetric part of the non-symmetric Laplacian is the Laplacian for a weighted undirected
graph and hence is an M-matrix \cite{BermanPlemmons} which is positive semidefinite.
Consider the weighted undirected graph with
adjacency matrix $\widetilde {\mathbf A} = (\bbg \Pi \mathbf P + \mathbf P^T \bbg \Pi)/2$.  
This is a weighted undirected graph with the same nodes as the original graph and an edge
whereever the original graph has an edge in either direction.
The vector of stationary probabilities for this graph is $\bbg \pi$,
proportional to the weighted degrees of the nodes in the new graph. 
The associated unnormalized Laplacian is
$\frc12(\mathbf L^{\mathrm r} +  (\mathbf L^{\mathrm r} )^T) $, which is therefore 
symmetric positive semi-definite with nullspace equal to $\mathsf{span}(\one )$
\cite{Chung97}.
The associated diagonally scaled Laplacian is
$\frc12(\mathbf L^{\mathrm d} +  (\mathbf L^{\mathrm d} )^T)
= \frc12
\bbg\Pi^{-\frc12} 
(\mathbf L^{\mathrm r} +  (\mathbf L^{\mathrm r} )^T)
\bbg\Pi^{-\frc12} 
$, which is therefore also symmetric positive semi-definite with
nullspace equal to  $\mathsf{span}(\sqrt{\bbg \pi})$.
The probability transition matrix for the new graph is
$\widetilde{\mathbf P} = \bbg \Pi^{-1} \widetilde {\mathbf A} = (\mathbf P + \bbg \Pi^{-1} \mathbf P^T \bbg \Pi)/2$.
Adding a symmetric rank-1 matrix ($\one  \one ^T$ or  $\sqrt{\bbg \pi}\sqrt{\bbg \pi}^T $, respectively)
makes the respective Laplacian matrices non-singular, moving the 0 eigenvalue to 
a positive number without moving the remaining eigenvalues. 
\QED

We can now notice that the Laplacian matrices $S(\mathbf L^{\mathrm r}) + \sqrt{\bbg \pi}\sqrt{\bbg \pi}^T $,
$S(\mathbf L^{\mathrm d}) + \one  \one ^T $ have just the right scaling 
to belong to a class of matrices for which GMRES (or any similar Krylov space minimum
residual method) has a guaranteed convergence rate.
We have the following theorem that is an immediate consequence of
Theorem \ref{GMRESconv} in the Appendix.
\boldparan{thm}{Corollary} \label{restartedthm}
Let $A$ be a real matrix such that $S(A)=(A+A^T)/2$ is symmetric positive definite and let
$\lambda_{\min}[S(A)]>0$ denote the smallest eigenvalue for $S(A)$.
The residual $\mathbf r_k$ obtained by restarted GMRES \cite{Saad-Iter} (restarting after $\ell$ inner steps)
after $k$ outer steps 
satisfies
$${
\frac{\|\mathbf r_k\|_2}{\|\mathbf r_0\|_2}
\le \left( 1 - \frac{(\lambda_{\min}[S(A)])^2} {\|A\|_2^2} \right)^{k\ell/2} 
}\myeqno{restartedbound}$$
{%
The cost of one outer step of restarted GMRES is
$\sfsl{cost}_1^{\rms{GMRES}} 
\ell\cdot\CMV + c_7(n\ell^2+\ell^{2})$.  }
{\boldpara{Proof}
According to Theorem \ref{GMRESconv}
the residual after $\ell$ steps of ordinary GMRES satisfies
$${
\frac{\|\mathbf r_{\ell}\|_2}{\|\mathbf r_0\|_2}
\le \left( 1 - \frac{(\lambda_{\min}[S(A)])^2} {\|A\|_2^2} \right)^{\ell/2} 
}\myeqno{gmresres}$$
Each time GMRES is restarted after $\ell$ steps, the residual is reduced
by the  factor in equ.\ (\ref{gmresres}).  After $k$ such repeats, the residual has been
reduced by a factor of at least (\ref{restartedbound}).
The cost estimate is based on an analysis of Algorithm \ref{alggmres},
detailed in \cite{Saad-Iter}, as sketched in the Appendix.  
}
\hfill\QED 

In summary, the total cost to find the vector of stationary probabilities is
$${
\apile[lll]{
\sfsl{cost}^{\rms{subspace}} &=& c^{\mathrm s} ((\ell + 1)\cdot\CMV + c_8(2n\ell^2+\ell^3)) \hfill
\\ &=& \frac{|\log(\sfsls{tol})|+|\log \|\mathbf r_0\| |}{|\log \rho|} ((\ell+1)\cdot\CMV + c_8(n\ell^2+\ell^3)),
}
}$$
where $\rho = {\lambda_{\ell+1}(\mathbf P)}$.
The cost to find each column of the pseudo-inverse with
a residual error of $\sfsl{tol}$, {given the vector of stationary probabilities}  is
$${
\apile[lll]{
\sfsl{cost}^{\rms{GMRES}} &=& 
c^{\mathrm g} (\ell\cdot\CMV + c_7(n\ell^2+\ell^{2}))
\\ &=&
\frac{|\log(\rms{tol})|+|\log \|\mathbf r_0\| |}{|\log \sigma|} (\ell\cdot\CMV + c_7(n\ell^2+\ell^{2}))
}}$$
where $\sigma = \left( 1 - \frac{(\lambda_{\min}[S(A)])^2} {\|A\|_2^2}
\right)^{\ell/2}$ {and 
$c^{\mathrm g} = \frac{|\log(\rms{tol})|+|\log \|\mathbf r_0\| |}{|\log \sigma|}$.}
{
Hence the total cost to obtain one column of the pseudo-inverse of an
Eulerian Laplacian equal to:
$${
\apile[lll]{
\sfsl{cost}^{\rms{one}} &=&{(c^{\mathrm s}\mbox{$+$}c^{\mathrm g}) } (\ell\cdot\CMV + c_9(n + n\ell^2+\ell^3))
,
}
}\myeqno{onecolcost}$$
where $\CMV = O(m)$ is proportional the number of edges in the graph.
{The space required is 
$\sfsl{space} =  O(\sfsl{\#edges} + n) +  (n \ell) + \ell^2)$,
where $\ell$ is the number of inner iterations of GMRES or the number of vectors
used in the subspace iteration (whichever is bigger).
}
We remark that in order to obtain the hitting time $h(i,k)$ for a given node
$k$, or to obtain the trust or influence measure
\cite{golnari2014revisiting,Liu16} for a given node $j$, only one column of
the pseudo-inverse is required at a cost shown in (\ref{onecolcost}).
Obtaining the entire pseudo-inverse, in cases where that would be
required, requires the computation of the stationary probabilities
only once, so the total cost would be bounded by:}
$${
\apile[lll]{
\sfsl{cost}^{\rms{total}} &\le& {(c^{\mathrm s} + n c^{\mathrm g})} (\ell\cdot\CMV + c_9(n + n\ell^2+\ell^3)) 
 . }
}$$

\nboldpara{par}{General Pseudo-Inverses}
\optpar
Here we show how to apply the previous to obtain the inverses or pseudo-inverses of 
general Laplacian matrices derived from strongly connected directed graphs. 
The approach
is to apply a row/column diagonal scaling to the non-Eulerian Laplacian to obtain a related Eulerian matrix
(similar to \cite{Cohen16}),
compute the pseudo-inverse using the previous methods, and unwind the diagonal scaling.
This can be applied to any Laplacian
$\widetilde {\mathbf L}$ that has all the following properties:
{\leftmargini=10ex
\begin{itemize}
\item[{\sf (Pa)}] $\widetilde {\mathbf L}$ is irreducible,
\item[{\sf (Pb)}] all the diagonal entries are strictly positive and all the off-diagonal entries are non-positive,
\item[{\sf (Pc)}] there is a strictly positive vector $\mathbf x$ so that $\widetilde {\mathbf L} \mathbf x = 0$.
\end{itemize} }
\noindent
Matrices satisfying property {\sf (Pb)} are Z-matrices {\cite{BermanPlemmons}}.
Alternatively, we can start with
an $(n-1)\times(n-1)$ matrix $\widetilde {\mathbf L}_{11}$ which shares properties {\sf (Pa)}
and {\sf (Pb)}, but has the property
{\leftmargini=10ex
\begin{itemize}
\item[{\sf (Pc')}] there is a strictly positive vector $\mathbf w$ so that $\widetilde {\mathbf L}_{11} \mathbf w $ is strictly positive.
\end{itemize} }%
\noindent
Then we use {Theorem} \ref{non-euler}(c) below to embed $\widetilde {\mathbf L}_{11}$ inside
an $n \times n$ matrix $\widetilde {\mathbf L}$ enjoying properties {\sf (Pa)}, {\sf (Pb)},
{\sf (Pc)}
and hence apply the same procedures.  Matrices satisfying {\sf (Pa)}, {\sf (Pb)}, {\sf (Pc')} are non-singular M-matrices and include
matrices that are strictly row-diagonally dominant with non-positive
off-diagonal entries \cite{BermanPlemmons}.
There are many other ways to characterize M-matrices (see \cite{BermanPlemmons}).

The pseudo-inverse of a diagonally scaled matrix is not the diagonally scaled pseudo-inverse of the original,
but the ordinary inverse of a diagonally scaled matrix is the diagonally scaled ordinary inverse of the original.
Hence one can apply the diagonal scaling to the leading principal submatrix of a Laplacian to map the problem
to the Eulerian scaling.
The following {theorem} provides a way to map from a matrix pseudo-inverse to the ordinary inverse of its principal submatrix
and vice versa, using only fast rank-one updates.
Alternatively, one can use $\bbg \Pi$ as a preconditioner on the unscaled Laplacian.

\boldparan{thm}{{Theorem}} \label{non-euler}
Let $C$ be an $n \times n$ non-singular matrix and suppose
$ A = C - \alpha \mathbf u \mathbf v\prm$
is singular with $A\mathbf u = A^T \mathbf v = \mathbf 0$, $\mathbf v^T\mathbf u = 1$
and $u_n > 0$, $v_n > 0$. 
Then
\begin{itemize}
\item[(a)]
the left and right nullspaces of $A$ are 
$\mathsf{nullsp}(A)=\mathsf{span}(\mathbf u)$, and
$\mathsf{nullsp}(A^T)=\mathsf{span}(\mathbf v)$.
The Moore-Penrose pseudo-inverse of $A$ is given as:\\

$${
\apile[l     l     l]{
A\psu &=& B ~ \defeq ~
  {\left(I_n - \frac1{\mathbf u^T\mathbf u}{\mathbf u \mathbf u^T} \right)}
C^{-1}
  {\left(I_n - \frac1{\mathbf v^T\mathbf v}{\mathbf v \mathbf v^T} \right)}
\\[1ex] &=&
C^{-1}
{- \frac{1}{\mathbf u^T\mathbf u} \mathbf u {\mathbf y^T }}
-\frac{       1    }{\mathbf v^T\mathbf v} \mathbf x \mathbf v^T
+\frac{\mathbf u^T        \mathbf x}{(\mathbf u^T\mathbf u)(\mathbf v^T\mathbf v)} \mathbf u \mathbf v^T ,
}}$$
\vspace*{-.5ex}
where
~~~ $\mathbf y^T ~=~ \mathbf u^T C^{-1}, ~~~ \mathbf x ~=~ C^{-1} \mathbf v$.
\item[(b)]
$A_{11}^{-1}$ exists and can be written in terms of $A\psu = B$, $\mathbf u$, $\mathbf v$:
$${
\apile[l     l     l]{
A_{11}^{-1} &=&
{\left(I_{n-1} + \frac1{u_n^2}{\mathbf u_1 \mathbf u_1^T}\right)}
B_{11}
{\left(I_{n-1} + \frac1{v_n^2}{\mathbf v_1 \mathbf v_1^T}\right)}
\\[1ex] &=& 
\left[ I_{n-1} ,~ -\mathbf u_1/u_n \right]
B
\left[\apile{ I_{n-1} \\ -\mathbf v_1^T/v_n} \right],
\\[2ex] &=&
   {B_{11} - \frac1{u_n}{\mathbf u_1 } \mathbf b_{21}^T - \frac1{v_n}\mathbf b_{12} {\mathbf v_1^T}
    + \frac{b_{nn}}{u_n v_n} \mathbf u_1 \mathbf v_1^T  }
}
}$$
where we have partitioned $B$ as above.
\item[(c)] We can write $A$ and $B = A\psu$ in terms of $A_{11}$, $\mathbf u$, $\mathbf v$
as follows
$${
\apile[l     l     l     l     l]{
A &=& 
{\left[ \apile[cc]{A_{11}  & -\frac1{u_n}A_{11} {\mathbf u_1}
       \\ -\frac1{v_n}{\mathbf v_1^T} A_{11} & \frac{\mathbf v_1^T A_{11} \mathbf u_1}{u_n v_n}
       }\right]}
&=& 
{\left[\apile{ I_{n-1} \\ -\frac1{v_n}{\mathbf v_1^T} } \right]}
A_{11}
{\left[ I_{n-1} ,~ -\frac1{u_n}{\mathbf u_1} \right]}
}}$$
\def\5{\rule{0pt}{0ex}}
$${\apile[l     l     l     l     l]{
B
&=&
  \left[\apile[cc]{
  B_{11} & \mathbf b_{12}
  \\
   \mathbf b_{21}^T
   & 
   b_{nn} 
  } \right]
&=& 
{\left[\apile{ I_{n-1} \mbox{$-$} \frac{\5\mathbf u_1\mathbf u_1^T}{\5\mathbf u^T\mathbf u}
 \\[.5ex] -\frac{u_n}{\mathbf u^T\mathbf u}\mbox{$\cdot$}\mathbf u_1^T} \right]}
A_{11}^{-1}
\left[ I_{n-1} \mbox{$-$} \frac{\5\mathbf v_1\mathbf v_1^T}{\mathbf v^T\mathbf v} ,
    -{\frac{v_n}{\mathbf v^T\mathbf v}}\mbox{$\cdot$}\mathbf v_1 \right] 
}
}$$
where the individual blocks of $B$ are 
$${
\apile[llll]{
\mbox{(a)}&B_{11} &=& A_{11}^{-1} - \mathbf u_1 \mathbf t^T 
- \mathbf w \mathbf v_1^T  + \frac{\mathbf u_1^T \mathbf w}{{\mathbf u^T}\mathbf u}\mbox{$\cdot$}\mathbf u_1 \mathbf v_1^T 
\\[1ex] \mbox{(b)} &
{\mathbf b_{12}} &=& {v_n\frac{\mathbf u_1^T \mathbf w }{\mathbf u^T\mathbf u} \mbox{$\cdot$}\mathbf u_1  - {v_n}\mathbf w}
\\[1ex]\mbox{(c)} &
\mathbf b_{21}^T &=&
   {u_n}{\frac{\mathbf u_1^T \mathbf w }{\mathbf u^T\mathbf u}} \mbox{$\cdot$}\mathbf v_1^T - {u_n}\mathbf t^T
\\[1ex] \mbox{(d)} &
b_{nn} &=& u_n {v_n} {\frac{\mathbf u_1^T \mathbf w }{\mathbf u^T\mathbf u}}
\\[1ex] \mc{4}{@{}l}{\mbox{where}
~~~ \mathbf w ~=~ \frac{{1}}{{\mathbf v^T\mathbf v}}A_{11}^{-1} \mathbf v_1 ,
~~~ \mathbf t^T ~=~ \frac{{1}}{{\mathbf u^T\mathbf u}}\mathbf u_1^T  A_{11}^{-1} }
}
}\myeqno{blocks}$$
\end{itemize}
{Part (a) is a special case of general theory of \cite{Meyer73},
and part (b) appears in \cite{boley11}.  Part (c) give a formula for the
Moore-Penrose pseudo-inverse that is similar to the formula in \cite{Meyer75}
for the so-called group inverse, but these two inverses agree only when the
left and right nullspaces match \cite{Roberts-GroupInv}.}
\boldpara{Proof}
The proof follows the same lines as that of {Theorem} \ref{psu}, after noting that
$AB = I - \frac{\mathbf v \mathbf v^T}{\mathbf v^T \mathbf v} $,
and $BA = I - \frac{\mathbf u \mathbf u^T}{\mathbf u^T \mathbf u}$.
\QED

Using this {theorem}, we briefly outline a feasible procedure to obtain the [pseudo]-inverse for an admissible
Laplacian matrix consisting of a sequence of diagonal scalings and the Eulerian Laplacian Algorithm \ref{alg1}. 

Suppose we have a matrix $\widetilde {\mathbf L} $ satisfying properties
{\sf (Pa)}, {\sf (Pb)}, {\sf (Pc)}, together with
a strictly positive vector $\mathbf x$ such that $0 = \widetilde {\mathbf L}  \mathbf x$.
{For example, an adjacency matrix $\widetilde {\mathbf A}$ for a strongly connected digraph,
with associated vector of out-degrees
$\tilde {\mathbf d}$, leads to the unnormalized Laplacian,
$\widetilde {\mathbf L} = \mathsf{Diag}(\tilde {\mathbf d}) - \widetilde {\mathbf A}$,
with $\mathbf x = \one$.}
Define
${
\widehat{\mathbf L} 
\defeq (\widehat {\mathbf D} - \widehat{\mathbf A} )
=  \widetilde{\mathbf L} \cdot \mathsf{Diag}(\mathbf x)
}$,
where $\widehat {\mathbf D}$ is the diagonal part of $\widehat{\mathbf L}$, and $-  \widehat{\mathbf A} $
has the rest.
This matrix
has the property that $\widehat {\mathbf L} \one  = 0$ while sharing the same 
left annihilating vector with $\widetilde{\mathbf L}$.
In this case, $\widehat {\mathbf P} = \widehat {\mathbf D}^{-1} \widehat {\mathbf A}$ is the probability
transition matrix for a random walk over this digraph, with {strictly positive} stationary probabilities $\bbg \pi$.
We can then {apply the diagonal scaling to $\widehat {\mathbf L}$ to obtain
the Eulerian Laplacian $\mathbf L^{\mathrm d} =
\bbg \Pi^{\frc12} (\widehat{\mathbf D}^{-1}\widehat{\mathbf L})\bbg \Pi^{-\frc12}$,
and then use Algorithm \ref{psu} to compute its pseudo-inverse.
The following algorithm uses this preprocessing to obtain the pseudo-inverse of
the original matrix $\widetilde {\mathbf L}$.}

{\def\theenumii{\roman{enumii}}
\boldparan{alg}{\VSP Algorithm} \label{alg2} \\
{\bf Input:} $\widetilde {\mathbf L} $ satisfying {\sf (Pa)}--{\sf (Pc)} and a strictly positive vector $\mathbf x$ such that
$\widetilde {\mathbf L} \mathbf x = 0$.\\
{\bf Output;} Pseudo-inverse of $\widetilde {\mathbf L}$.

\begin{enumerate}
\itemsep=0pt \topsep=0pt 

\item Form ${ \widehat{\mathbf L}
           \defeq
    \widetilde{\mathbf L} \cdot \mathsf{Diag}(\mathbf x) }
    =
	   (\widehat {\mathbf D} - \widehat{\mathbf A} )$,
    where $\widehat {\mathbf D} = \widetilde{\mathbf D}\cdot\mathsf{Diag}(\mathbf x)$ is the diagonal part of $\widehat {\mathbf L}$,
    and $\mathbf A$ has the rest.

\item  Form the probability transition matrix 
   $\widehat {\mathbf P}^T = (\widehat {\mathbf D}^{-1}  \widehat{\mathbf A} )^T$.

\item (*)  Use Algorithm \ref{alg1} to compute the vector of stationary probabilities $\bbg \pi$ and
   the pseudo-inverse $(\mathbf L^{\mathrm d})\psu$ 
   of  $\mathbf L^{\mathrm d} ~\defeq~ \bbg \Pi^{\frc12} (I \mbox{$-$} \widehat{\mathbf D}^{-1}\widehat{\mathbf A})\bbg \Pi^{-\frc12} 
\mbox{$=$} \bbg \Pi^{\frc12} (\widehat{\mathbf D}^{-1}\widehat{\mathbf L})\bbg \Pi^{-\frc12}$.

\item (*) Use Theorem \ref{psu}(b) to obtain the $(n-1)\times(n-1)$ matrix
$(\mathbf L_{11}^{\mathrm d})^{-1}$ from $(\mathbf L^{\mathrm d})\psu$.
Here $\mathbf u = \sqrt{\bbg \pi}$.

\item (*) Form the $(n-1) \times (n-1) $ matrix
$(\widetilde {\mathbf L}_{11})^{-1} = 
(\mathsf{Diag}(\mathbf x_{1:n-1}))
\bbg \Pi_1^{-\frc12}
(\mathbf L_{11}^{\mathrm d})^{-1} 
\bbg \Pi_1^{\frc12}
\widehat{\mathbf D}_1^{-1}
$.
\item \label{a2last} (*) Use Theorem \ref{non-euler}(c) to obtain $\widetilde {\mathbf L}\psu $ from 
$(\widetilde {\mathbf L}_{11})^{-1}$:  
Use (\ref{blocks}) with
$A_{11}^{-1}=(\widetilde {\mathbf L}_{11})^{-1}$,
$\mathbf u = \mathbf x $, and $\mathbf v = \widetilde{\mathbf D}^{-1}\bbg \pi$,
i.e., $v_i = \pi_i/(\tilde d_i)$, $i=1,\ldots,n$.
The vectors $\mathbf u ,\mathbf v$ are the right and left annihilating vectors for the original
Laplacian $\widetilde {\mathbf L}$.

\end{enumerate}
}%
{%
\noindent
In the steps marked (*), if only a limited set of columns are required, only those columns must be computed,
though the entire vector $\bbg \pi$ must be computed.
For example, if only column $j \le n-1$ is required,
then only column $j$ must be computed in the steps marked (*),
and we do not need the entire vector
$\mathbf t$, but only its $j$-th entry 
$t_j = \frac1{\mathbf u^T\mathbf u} \mathbf u_1^T [A_{11}^{-1}]_{:j}$, obtainable from the $j$-th column of
$A_{11}^{-1} = (\widetilde {\mathbf L}_{11})^{-1}$ already computed.
In any case step \ref{a2last} requires the solution of an extra system of
linear
equations for $\mathbf w$ (\ref{blocks}):
$${\mathbf w = \frac1{\mathbf v^T \mathbf v} (\widetilde{\mathbf L}_{11})^{-1} \mathbf v_1
= 
\frac1{\mathbf v^T \mathbf v}
(\mathsf{Diag}(\mathbf x_{1:n-1}))
\bbg \Pi_1^{-\frc12}
\fbox{$(\mathbf L_{11}^{\mathrm d})^{-1} 
\mathbf z_1$}}$$
where
${
\mathbf z_1 \defeq
\bbg \Pi_1^{\frc12}
\widehat{\mathbf D}_1^{-1}
\mathbf v_1 
}$.
Using (\ref{AA}), the boxed expression can be computed as follows:
$${
(\mathbf L_{11}^{\mathrm d})^{-1} \mathbf z_1 
=
\left[ I_{n-1} ,~ -\sqrt{\bbg \pi_1}/\sqrt{\pi_n} \right]
(\mathbf L^{\mathrm d})\psu
\left[\apile{ I_{n-1} \\ -\sqrt{\bbg \pi_1}^T/\sqrt{\pi_n}} \right]
\mathbf z_1
\defeq
\left[ I_{n-1} ,~ -\sqrt{\bbg \pi_1}/\sqrt{\pi_n} \right]
\fbox{$(\mathbf L^{\mathrm d})\psu \mathbf z$}
}$$
where $\mathbf z^T \defeq \left[\mathbf z_1^T  , -\frac{\sqrt{\bbgs{\pi_1}}^T\mathbf z_1}{\sqrt{\pi_n}}\right] $.
The expression 
$(\mathbf L^{\mathrm d})\psu \mathbf z$
can be computed as in Algorithm \ref{alg1} step \ref{a1step5} using the same 
restarted GMRES procedure, based on the identity (\ref{LapInv}).
}

\nboldpara{par}{Performance}
\optpar
To illustrate how the theoretical complexity corresponds to practice
we generate a sequence of synthetic graphs using preferential attachment
\cite{AlbertBarabasi02} with $2n$ edges plus an extra set of $n$
randomly placed one-way edges to make the graph a digraph
{for a total of $3n$ edges. }
Table \ref{table1} shows 
the time to compute the stationary probabilities $\bbg \pi$
and the time to solve a single linear system
involving $\mathbf L^{\mathrm d}$ using restarted GMRES.
{We run the methods in their original unpreconditioned form to show
the correspondence between the theoretical complexity bounds and the
numerical behavior in practice. }
We also show the number of matrix-vector products, which is solely a function
of the number of overall iterations, which in turn is entirely dependent
on the eigenvalue distribution of Laplacian.  This, of course, depends on
the nature of the underlying graph and would have to be analysed on a
case-by-case basis.  In the sequence of synthetic graphs constructed for this
illustration, it is seen that the number of iterations is a slowly growing
function of the dimension.  Except for the modest increase in the
number of matrix-vector products, the cost of the methods
approximately double when the dimension $n$ is doubled. 
Here, the iterative methods were applied with a zero tolerance of $10^{-9}$.
Using Matlab R2018a on a 2.5GHz Linux desktop with 8 Intel(R) Core(TM) i7-7700 CPU cores and 32 GB RAM,
each experiment was repeated 4 times with averages shown in Table \ref{table1}.  

{The computation requires the storage
of the adjacency matrix in sparse format plus
up to $\ell$ auxiliary vectors of length $n$ and an $\ell \times \ell$
matrix, where $\ell$ is a user parameter
independent of $n$. }
{The off-the-shelf methods to compute the pseudo-inverse would use direct methods: solving
the homogeneous system (\ref{directeig}) for the vector of stationary probabilities, and using
Theorem \ref{psu}(c) or \ref{non-euler}(c) to solve for the pseudo-inverse given the left and right
annihilating vectors $\mathbf u,\mathbf v$.  In both cases, the process would involve solving a non-symmetric
system of linear equations using a method like Gaussian elimination with partial pivoting.
This would require $O(n^3)$ work with $O(n^2)$ storage, even if solving for just one column.
For example, the largest case shown ($n=2^{18}$) in the table would require $n^2 \cdot 8$ $\approx$ 550GB space,
while the sparse iterative methods proposed here required only 
$8$ $\cdot$ ($n\cdot {\mathrm d} \cdot n\ell$ $+$ l.o.t)
$\approx$ 71MB, where we use ${\mathrm d} = 4 = 1$ $+$ \#edges\_per\_vertex,
$\ell = 30$, and 8 bytes per double word.
Results using purely direct methods were reported in \cite{BSMDMA2019}, where an
off-the-shelf minimum degree ordering was
used to reduce the fill in exact Gaussian elimination,
but the induced
fill was still observed to be $O(n^2)$ for the randomaly generated synthetic graphs.
One could trade off accuracy for the direct methods by fixing the ordering of the equations to reduce
fill, thereby substantially reducing the cost as proposed in \cite{Cohen18,Cohen16}.
However, here we avoid the issue of fill-in entirely by using purely iterative methods
with guarantees on the progress made at
each iteration toward the solution.   One can iterate to reduce the error to 
any desired tolerance within the range of the underlying arithmetic precision,
without adding to the memory footprint.
}

Table \ref{table1} also shows performance on an Epinions
data set \cite{Epinions}, augmented with an evaporating node
with an evaporating probability of $\gamma=0.05$ followed by a uniformly random
restart
(akin to the pagerank teleportation probability \cite{PageBrin98}).
In other words, at every
transition the walker has a 5\% chance of ``evaporating'' to the extra node, 
and from the extra node the walker transitions to one of the original nodes
with equal probability. Including the extra node and associated links and 9 old
nodes not otherwise connected to any other old node, the resulting graph has 75,889 nodes and 660,613 links.

\smallskip

\begin{table}
	\caption{Performance on synthetic graphs \& one social network } \label{table1}
\centering 
\begin{tabular}{|r|r|r|r|r|} 
 \hline
  \mc{1}{|c}{$n$}  
	& \mc{2}{|c}{get \(\bbg\pi\)}  & \mc{2}{|c|}{GMRES (1 col)}
\\ \hline
\mc{1}{|c}{dim}
	& \mc{1}{|c}{\sfsl{\#Mv}}
	& \mc{1}{|c|}{ time  (ms)}
	& \mc{1}{|c}{\sfsl{\#Mv}}
	& \mc{1}{|c|}{ time  (ms)}
\\ \hline
     1024 & 237 &     3.919       & 59  &     7.097         \\
     2048 & 286 &     5.297       & 65  &     9.331         \\
     4096 & 303 &    11.685       & 68  &    23.451         \\
     8192 & 369 &    28.223       & 74  &    51.270         \\
    16384 & 429 &    56.804       & 82  &    82.047         \\
    32768 & 347 &    88.941       & 77  &   135.895         \\
    65536 & 391 &   192.741       & 83  &   265.114         \\
   131072 & 429 &   427.461       & 86  &   483.605          \\
   262144 & 528 &  1130.698       & 95  &  1136.207         
\\ \hline \hline
  \mc{5}{|c|}{Epinions data set}
\\ \hline
   75889  & 682 &   477.362       & 68  &   304.657        
\\ \hline
\end{tabular} %
\end{table}

\nboldpara{par}{Discussion and Conclusions}
\optpar
We have used several off-the-shelf matrix iterative methods to compute 
individual columns of the pseudo-inverse of a digraph Laplacian matrix in time linear in the
number of edges times 
a factor depending on how difficult it is
to cut the graph into separate connected components. 
The full pseudo-inverse can be computed in time that is amortized to
almost constant time per matrix entry.  
The methods proposed here depend exclusively on iterative methods and do not make any use
of an elimination scheme that results in fill-ins, unlike methods using a variant of Gaussian
elimination.  They are relatively efficient and enjoy a plethora of available implementations.

The overall complexity is entirely dependent on the convergence rates for the
iterative
eigensolver and GMRES, which in turn
depend on
the smallest nonzero eigenvalues of
$\mathbf I - \mathbf P$ and $\mathbf L+\mathbf L^T$, respectively.
The former is related to the mixing rate of the random walk, while the latter
is closely related to the Cheeger constant \cite{ChungDigraph}.
In both cases a tiny eigenvalue corresponds 
to a graph that can be split with a small cut.
The convergence of any iterative method would depend on similar quantities in some fashion.

\section*{Appendix}

We collect in this Appendix some results from the literature on which this paper is based.

\nboldpara{par}{Compute Stationary Probabilities}
The vector of stationary probabilities $\bbg \pi$ is the eigenvector of $\mathbf P^T$ corresponding to the
eigenvalue $\lambda=1$.  Since the underlying graph is strongly connected, the Perron Frobenius theory
guarantees eigenvalue $\lambda=1$ is simple.  The number of other eigenvalues of modulus 1 is equal to the
periodicity of the graph or random walk.  For instance, a bipartite graph will have an eigenvalue
$-1$.  If $\sfsl{per}$ is the period of the graph and we use $\ell>\sfsl{per}$ vectors in the following algorithm
then the algorithm is guaranteed to converge at a rate bounded by $|\lambda_{\ell+1}(P)| < 1$
\cite{Stewart} since 
$\lambda=1$ is known and is a simple eigenvalue of largest modulus.

\boldparan{alg}{\VSP Algorithm} \label{algsub} {\bf Modified Subspace Iteration.} \cite{Stewart}
\\ {\bf Input:} matrix $A$, hyperparameters: \textsf{tol}, 
initial guess $X^{[0]} \in \mathbb{R}^{n \times \ell}$.
\\ {\bf Output:} $\mathbf z = $ eigenvalue corresponding to eigenvalue 1.
{\leftmarginii=4ex
\begin{shortenum}
  \item Set $Z = \mbox{orthogonalize}(X^{[0]})$, where $Z_{:,1}$ is all non-negative.
  \item Repeat until convergence: 
  \begin{shortenum}
    \item Set $Q = \mbox{orthogonalize}(A Z)$.
    \item Compute $\ell \times \ell$ Schur Decomposition $[U T U^T] = Q^T A Q$
    with the diagonal entries of $T$ ordered to put the entry closest to 1 in the 1,1 position.
    \item Set $Z=Q U$.  Ensure the first column $Z_{:,1}$ is all non-negative
    (flipping signs of rows of $Z$ to make the first column all non-negative, if necessary).
  \end{shortenum}
  \item Return $\mathbf z = Z_{:,1}$.
\end{shortenum}}
{
The cost per iteration is $\CMV\cdot(\ell+1) + c_5(n\ell^2) + c_6(\ell^3)$ where
the first term accounts for the matrix vector products,
the second term
accounts for the orthogonalization (Alg \ref{algsub} step 2(i)) and the third term accounts for the $\ell\times\ell$ Schur decomposition (step 2(ii)).
Here 
$\sfsl{Mv}$ is the number of matrix vector products.  Each matrix-vector product requires
$\CMV = 2 \cdot \mathsf{nnz}(\sfsl{matrix})$ flops (one multiply and one add per matrix entry).
Here each matrix entry corresponds exactly to an edge in the graph.
So the total cost per iteration is
$${
\sfsl{cost}_1^{\rms{subspace}} =
\underbrace{c_2 \sfsl{\#nonzeros} \cdot ( \ell+1)}_{{\CMV}} + \underbrace{c_3 n
\ell^2}_{{\textrm{orthogonalization}}} + \underbrace{c_4 \ell^3}_{{\textrm{Schur decomp}}},
}$$
for some small constants $c_2,c_3,c_4$ at most 10.
The storage required is (in words)
$${
\sfsl{space}^{\rms{subspace}} =  \underbrace{(\sfsl{\#nonzeros} + n)}_{\textrm{sparse }A} + \underbrace{ 2 (n
\ell)}_{\textrm{iterates}Z,Q} + \underbrace{O(\ell^2)}_{\textrm{temporaries }U,T}
}\myeqno{spacesub}$$}

\nboldpara{par}{Restarted GMRES}
The heart of the computation of the pseudo-inverse is the use of {Theorem} \ref{psu}(a)
to convert a pseudo-inverse computation to an ordinary inverse computation.  
The restarted GMRES algorithm has received much attention in the literature (see \cite{Saad-Iter} and references
therein) with many enhancements for numerical stability that do not impact the cost by more than a constant
factor.  For the purposes of showing the overall cost of the algorithm, we show a simplified sketch of the
basic algorithm.  By using restarted GMRES, as opposed to ordinary GMRES, we bound the cost of
each iteration and the memory footprint.

\boldparan{alg}{\VSP Algorithm} \label{alggmres} {\bf Arnoldi-based Restarted GMRES.}\\
{\bf Input:} Matrix $A$, right hand side $\mathbf b$,
hyperparameters: restart count $\ell$, outer iteration limit \textsf{maxit}, tolerance
\textsf{tol}, initial vector $\mathbf x^{[0]}$.\\
{\bf Output:} solution $\mathbf x$ such that $\|\mathbf r\| = \|A \mathbf x - \mathbf b\| < \sfsl{tol}$.
{\leftmarginiii=5ex
\leftmarginiv=5ex
\begin{shortenum}
  \item Compute $\mathbf r^{[0]} = \mathbf b - A \mathbf x^{[0]}$
  \item For $k = 0,1,2,\ldots \sfsl{maxit}$: \label{outer}
  \begin{shortenum}
    \item Set $\beta = \|\mathbf r^{[k]}\|_2$ and set $\mathbf v_1 = \mathbf r^{[k]}/\beta$.
    \item If $\beta < \sfsl{tol}$, return solution $\mathbf x = \mathbf x^{[k]}$.
    \item Generate orthonormal Arnoldi basis $V_{\ell+1} = [\mathbf v_1,\cdots,\mathbf v_{\ell+1}]$
    for the Krylov space \\\relax  $\mathsf{span}\{\mathbf r_0, A\mathbf r_0 ,\ldots,A^{l}\mathbf r_0\}$, and upper Hessenberg matrix
    $\bar H_\ell \in \mathbb{R}^{\ell+1 \times \ell}$ such that
    $A V_\ell = V_{\ell+1} \bar H_{\ell}$.
    \item Compute $\mathbf y^{[k]} = \mathrm{arg } \min_{\mathbf y} \|\beta \mathbf e_1 - \bar H_\ell \mathbf y\|_2^2$.
    \item Set $\mathbf x^{[k+1]} = \mathbf x^{[k]} + V_\ell \mathbf y^{[k]}$ 
  \end{shortenum}

\end{shortenum}
}

The cost of one outer step \ref{outer} of restarted GMRES is $\ell\cdot\CMV + O(n\ell^2+\ell^{2})$ \cite{Saad-Iter}.
Here $\CMV$ is the cost of one matrix vector product involving sparse matrix $A$.  This takes one
floating multiply and one floating add for each nonzero matrix element.  So the cost is
$\CMV = O(\mathsf{nnz}(A))$.
{The Arnoldi step 2(iii) has one matrix vector
product and an orthogonalization step for
each of the $\ell$ Krylov vectors generated \cite{Saad-Iter}.}
Step 2(iv) is an $(\ell+1) \times \ell$ least squares problem costing
$O(\ell^2)$ to solve, due to the special Hessenberg structure of $\bar H$.
The number of outer iterations required is controlled by the eigen-structure
of the symmetric part $S(A)=(A+A^T)/2$, which is related to the separability of the underlying graph
\cite{ChungDigraph}.  
{Hence the total cost of restarted GMRES is
$${
\apile[lll]{
\sfsl{cost}^{\rms{GMRES}} &=& 
c^{\mathrm g} (\ell\cdot\CMV + c_7(n\ell^2+\ell^{2}))
\\ &=&
\underbrace{\frac{|\log(\rms{tol})|+|\log \|\mathbf r_0\| |}{|\log \sigma|}}_{{\textrm{\#iterations}}}
(\underbrace{\ell\cdot\CMV}_{{\textrm{Mat\raisebox{-.5ex}{*}vec}}} +
c_7(\underbrace{n\ell^2}_{{\textrm{orthogonalization}}}
+\underbrace{\ell^{2}}_{{\textrm{work with }H}}))
}}$$
where $\sigma = \left( 1 - \frac{(\lambda_{\min}[S(A)])^2} {\|A\|_2^2}
\right)^{\ell/2}$ {and 
$c^{\mathrm g} = \frac{|\log(\rms{tol})|+|\log \|\mathbf r_0\| |}{|\log \sigma|}$.}
The storage required is (in words)
$${
\sfsl{space}^{\rms{GMRES}} =  \underbrace{(\sfsl{\#nonzeros} + n)}_{\textrm{sparse }A} + \underbrace{ 2 (n
\ell)}_{\textrm{Krylov vectors }V} + \underbrace{O(\ell^2)}_{\textrm{temporaries }H,Y}.
}\myeqno{spacegmres}$$
}

Regarding the number of GMRES iterations,
we have the following bound which yields Corollary \ref{restartedthm} as an
immediate consequence.

\boldparan{thm}{Theorem} \label{GMRESconv} \cite{EE,Elman,Joubert94,liesen2012field},
Let $A$ be a matrix such that $S(A) \defeq (A+A^H)/2$ is Hermitian positive definite and let
$\lambda_{\min}[S(A)]>0$ denote the smallest eigenvalue for $S(A)$.
The residual $\mathbf r_k$ obtained by GMRES \cite{Saad+Schultz} after $k$ steps 
applied to the linear system $A \mathbf x = b$
satisfies
\vspace*{-.75ex}
$${
\frac{\|\mathbf r_k\|_2}{\|\mathbf r_0\|_2}
\le \left( 1 - \frac{(\lambda_{\min}[S(A)])^2} {\|A\|_2^2} \right)^{k/2} 
}$$
\hfill\raisebox{2ex}[0pt][0pt]{\QED}

{To be self-contained,}
we give a sketch of a proof for this essential bound, referring to
to \cite{EE,Elman,Joubert94,liesen2012field} for detailed proofs for
this and several tighter bounds.  First we need the following Lemma
\boldparan{thm}{Lemma} \label{Lemmap1}
Let $A \in \mathbb{C}^{n \times n}$ such that $S(A)=(A+A^H)/2$ is positive definite,
and let $\mathbf v \in \mathbb{C}^n$ with $\|\mathbf v\|_2 = 1$ be given.
Let $\alpha_{\mathbf v} = [(A \mathbf v)^H \mathbf v]/[(A \mathbf v)^H (A \mathbf v)]$.  Then
$${
\|\mathbf v - \alpha_{\mathbf v} A \mathbf v\|^2
\le 1 - \lambda_{\min}[S(A^{-1})] \cdot \lambda_{\min}[S(A)]
\le 1 -  \frac{(\lambda_{\min}[S(A)])^2}{\|A\|^2}
}\myeqno{lemmap1bound}$$
\boldpara{Proof}
$\alpha_{\mathbf v}$ is the value achieving the minimum in the scalar least squares problem
$\min_{\alpha}\|\mathbf v - \alpha A \mathbf v\|_2^2$ and hence satisfies the Galerkin condition
$(A \mathbf v)^H (\mathbf v - \alpha_{\mathbf v} A \mathbf v) = 0$.  So we have
$${
\|\mathbf v - \alpha A \mathbf v\|^2 =  \mathbf v^H (\mathbf v - \alpha_{\mathbf v} A \mathbf v) = 1 - \alpha_{\mathbf v} \cdot \mathbf v^H A \mathbf v 
= 1 - \frac{\mathbf w^H A^{-1} \mathbf w}{\mathbf w ^H \mathbf w} \cdot \mathbf v^H A \mathbf v ,
}$$
where $\mathbf w = A \mathbf v$.
A well known result on field of values for any matrix $M$ whose Hermitian part $S(M)=(M+M^H)/2$ is positive definite
is the inequality \cite{hojo:91} for any $\mathbf x \not= 0$
$${
\left|\frac{\mathbf x^H M \mathbf x}{\mathbf x ^H \mathbf x} \right | \ge 
Re \left(\frac{\mathbf x^H M \mathbf x}{\mathbf x ^H \mathbf x} \right ) 
\ge \lambda_{\min}[S(M)] > 0.
}$$
Hence the first inequality (\ref{lemmap1bound}) follows. 
The remaining inequality follows from the identity
$${ S(A^{-1}) = A^{-1} \cdot S(A) \cdot A^{-H} . }$$
Inverting both sides and taking norms yields
$${
\frac{\DS 1}{\DS \lambda_{\min}[S(A^{-1})]}
= \|[S(A^{-1})]^{-1}\|_2 \le \|A\|_2^2 \|[S(A)]^{-1}\|_2 = 
\frac{\DS \|A\|_2^2}{\DS \lambda_{\min}[S(A)]}
}$$
\hfill\QED

\boldpara{Sketch of Proof of Theorem \ref{GMRESconv}.}
We include a proof sketch stripped down to its bare essentials.
GMRES on a matrix $A$ with initial residual $\mathbf r_0$ will find in $k$ steps a solution
with a residual $\mathbf r_k$ satisfying 
$\|\mathbf r_k\| = \min_{p \in \mathbb{P}_k}\|p(A)\mathbf r_0\|$,
where $\mathbb{P}_k$ is the set of all polynomials $p$ of degree up to $k$ satisfying $p(0)=1$.
In particular, after a single step 
$\|\mathbf r_1\| = \min_{p \in \mathbb{P_1}}\|p(A)\mathbf r_0\| \le 
\|(I - \alpha(\mathbf r_0)) A \mathbf r_0\|$,  where $\alpha(\mathbf r_0) = 
[(A \mathbf r_0)^H \mathbf r_0]/[(A \mathbf r_0)^H (A \mathbf r_0)]$. 
Hence we have the bound from Lemma \ref{Lemmap1}:
$\|\mathbf r_1\|_2^2 \le \|\mathbf r_0\|_2^2 \cdot [1 - \lambda_{\min} S(A) \lambda_{\min} S(A^{-1})]$.
This amounts to a single step of a variant of the classical Richardson iteration.
Repeating this Richardson iteration yields
$${
\mathbf r_k^{\mathrm{[rich]}}
=
\left(I - \alpha( \mathbf r_{k-1}^{\mathrm{[rich]}}) A\right) 
\cdots
\left(I - \alpha( \mathbf r_{1}^{\mathrm{[rich]}}) A\right) 
\left(I - \alpha( \mathbf r_{0}^{\mbox{~}}) A\rule{0pt}{2.5ex}\right) \mathbf r_0 .
}$$
The norm of the residual after $k$ Richardson steps is
bounded above by the convergence bound and bounded below by the norm of the GMRES residual
after $k$ steps:
$${\|\mathbf r_k^{\mathrm{[GMRES]}}\|_2^2 =  \min_{p \in \mathbb{P}_k} \|p(A) \mathbf r_0\|_2^2
\le
\|\mathbf r_k^{\mathrm{[rich]}}\|_2^2
\le 
\left( 1 - \lambda_{\min}[S(A)]\lambda_{\min}[S(A^{-1})]\right)^{k}
}$$
\hfill\raisebox{2ex}[0pt][0pt]{\QED}

\section*{Acknowledgements}
This research was supported in part by NSF grants 1835530 and 1922512.  The author would like to thank
the anonymous reviewer for many helpful comments that greatly improved the manuscript.


\end{document}